\documentclass[preprint, review,12pt]{elsarticle}

\usepackage{amssymb}
\usepackage{amsmath}

\journal{Number Theory}
\begin{document}

\begin{frontmatter}

%\preprint{}

\title{On Apery's Constant and Catalan's Constant}% Force line breaks with \\

\author{Akhila Raman }

\address{University of California at Berkeley, CA-94720. Email: akhila.raman@berkeley.edu. Ph: 510-540-5544}

%\affiliation{University of California at Berkeley, CA-94720. Email: akhila.raman@berkeley.edu. Ph: 510-540-5544}


\begin{abstract}
In this paper, Riemann's Zeta function with odd positive integer argument is represented as an
infinite summation of integer powers of $\pi$ with rational coefficients. Specific
values for Apery's Constant and Catalan's Constant are then derived.

\end{abstract}

\begin{keyword}



\end{keyword}

\end{frontmatter}
   
%\maketitle


%\section{}
%\label{Introduction}



%% References

\section{\label{sec:level1}Introduction\protect\\  \lowercase{} }

It is well known that Riemann's Zeta function with even positive integer 
argument is given by 

\begin{equation} \zeta{(2n)} = \sum_{k=1}^{\infty} \frac{1}{k^{2n}} =  \frac{B_{2n}}{2}  (2\pi)^{2n}\frac{(-1)^{n+1}}{ !({2n})}  \end{equation}

where $B_{2n}$ is a Bernoulli number. No such simple expression is known for odd positive integer argument. In this section, Riemann's Zeta function with odd positive integer argument is represented as an
infinite summation of integer powers of $\pi$ with rational coefficients.
\\

It is well known from the theory of divergent series[1] that in the interval $|\theta|<\pi$,
\begin{equation} \sin(\theta)-\sin(2\theta)+ \sin(3\theta)-\sin(4\theta)+...=(1/2)\tan(\theta/2)  \end{equation}

The Right Hand Side of this equation can be expanded using Maclaurin's series[2] as follows:

\begin{equation}  (1/2)\tan(\theta/2) = (1/2) \sum_{n=1}^{\infty}{c_{n}(\theta/2)^{2n-1}} \end{equation}

where the coefficients $c_{n}$ are given in terms of Bernoulli numbers $B_{2n}$
as $c_{n}=B_{2n} (2^{2n}-1) 2^{2n}\frac{(-1)^{n+1}}{!({2n})}$.\\

%\begin{equation}   \end{equation}
$\textbf{{Step 1}}$   \\ 
Integrating Eq.2  over $ [0,  \theta]$ in Step 1 [Integrating over $ [0,  \phi]$ and substituting $\phi = \theta$ ],  %% REDO eqn

\begin{equation}   -[\frac{\cos(\theta)}{1}-\frac{\cos(2\theta)}{2}+ \frac{\cos(3\theta)}{3}-\frac{\cos(4\theta)}{4}+...]+[\frac{1}{1}-\frac{1}{2}+\frac{1}{3}-\frac{1}{4}+...] 
% \end{equation}
%\begin{equation} 
= (1/2) \sum_{n=1}^{\infty}{\frac{c_{n}(\theta)^{2n}}{(2n) 2^{2n-1}}  }  %% REDO EQN
\end{equation}

Putting  $ \theta = \frac{\pi}{2} $ in Eq.4, we  get,  \\
\begin{equation}  A_{1}= [\frac{1}{1}-\frac{1}{2}+\frac{1}{3}-\frac{1}{4}+...]  =  \sum_{n=1}^{\infty}{ E_{n}(1) ({\frac{\pi}{2}})^{2n}  }  \end{equation}\\

where  $  E_{n}(1) = D_{n}(1)$ and $D_{n}(k) = \frac{c_{n}}{2^{2n-1}(2n+k-1)P_{k}} $
and  $D_{n}(k) = \frac{D_{n}(k-1)}{(2n+k-1)} $.\\

\begin{equation}   [\frac{\cos(\theta)}{1}-\frac{\cos(2\theta)}{2}+ \frac{\cos(3\theta)}{3}-\frac{\cos(4\theta)}{4}+...] = A_{1} - \frac{1}{2} \sum_{n=1}^{\infty}{ D_{n}(1) (\theta)^{2n}  }   \end{equation}

$\textbf{{Step 2}}$   \\ 
Integrating Eq.6  over  $ [0,  \theta]$  in Step 2, 

\begin{equation}   [\frac{\sin(\theta)}{1^{2}}-\frac{\sin(2\theta)}{2^{2}}+ \frac{\sin(3\theta)}{3^{2}}-\frac{\sin(4\theta)}{4^{2}}+...] = A_{1}\theta - \frac{1}{2} \sum_{n=1}^{\infty}{ D_{n}(2) (\theta)^{2n+1}  }   \end{equation}

Putting  $ \theta = \frac{\pi}{2} $ in Eq.7, we  get,  \\
\begin{equation}  A_{2}= [\frac{1}{1^{2}}-\frac{1}{3^{2}}+\frac{1}{5^{2}}+...]  = A_{1} \frac{\pi}{2} - \frac{1}{2} \sum_{n=1}^{\infty}{ D_{n}(2) ({\frac{\pi}{2}})^{2n+1}  }  \end{equation}\\

\begin{equation}  A_{2}=  \sum_{n=1}^{\infty}{ E_{n}(2) ({\frac{\pi}{2}})^{2n+1}  }    \end{equation}
where  $  E_{n}(2) = E_{n}(1) - \frac{D_{n}(2)}{2}  $ \\

$\textbf{{Step 3}}$   \\ 
Integrating Eq.7  over  $ [0,  \theta]$  in Step 3, 

\begin{equation}   -[\frac{\cos(\theta)}{1^{3}}-\frac{\cos(2\theta)}{2^{3}}+ \frac{\cos(3\theta)}{3^{3}}-\frac{\cos(4\theta)}{4^{3}}+...]+[\frac{1}{1^{3}}-\frac{1}{2^{3}}+\frac{1}{3^{3}}-\frac{1}{4^{3}}+...] =  A_{1}\frac{\theta^{2}}{!{2}} - \frac{1}{2} \sum_{n=1}^{\infty}{ D_{n}(3) (\theta)^{2n+2}  }  
\end{equation}

Putting  $ \theta = \frac{\pi}{2} $, we  get,  \\
\begin{equation}  A_{3}= [\frac{1}{1^{3}}-\frac{1}{2^{3}}+\frac{1}{3^{3}}-\frac{1}{4^{3}}...]  = \frac{1}{1-\frac{1}{2^{3}}}[\frac{A_{1}}{!{2}} (\frac{\pi}{2})^{2} - \frac{1}{2} \sum_{n=1}^{\infty}{ D_{n}(3) ({\frac{\pi}{2}})^{2n+2}  } ] \end{equation}\\

\begin{equation}  A_{3}=  \sum_{n=1}^{\infty}{ E_{n}(3) ({\frac{\pi}{2}})^{2n+2}  }    \end{equation}
where  $  E_{n}(3) = \frac{ - \frac{D_{n}(3)}{2} + \frac{E_{n}(1)}{!{2}}}{1-\frac{1}{2^{3}}}$ \\

we have \\
\begin{equation}   [\frac{\cos(\theta)}{1^{3}}-\frac{\cos(2\theta)}{2^{3}}+ \frac{\cos(3\theta)}{3^{3}}-\frac{\cos(4\theta)}{4^{3}}+...] = A_{3} - \frac{A_{1} \theta^{2}}{!{2}}+ \frac{1}{2} \sum_{n=1}^{\infty}{ D_{n}(3) (\theta)^{2n+2}  }   \end{equation}

$\textbf{{Step 4}}$   \\ 
Integrating Eq.13  over  $ [0,  \theta]$  in Step 4, 

\begin{equation}   [\frac{\sin(\theta)}{1^{4}}-\frac{\sin(2\theta)}{2^{4}}+ \frac{\sin(3\theta)}{3^{4}}-\frac{\sin(4\theta)}{4^{4}}+...] = A_{3}\theta  - \frac{A_{1} \theta^{3}}{!{3}}+ \frac{1}{2} \sum_{n=1}^{\infty}{ D_{n}(4) (\theta)^{2n+3}  }   \end{equation}

Putting  $ \theta = \frac{\pi}{2} $, we  get,  \\
\begin{equation}  A_{4}= [\frac{1}{1^{4}}-\frac{1}{3^{4}}+\frac{1}{5^{4}}+...]  = A_{3} \frac{\pi}{2} - \frac{A_{1} (\frac{\pi}{2})^{3}}{!{3}} + \frac{1}{2} \sum_{n=1}^{\infty}{ D_{n}(4) ({\frac{\pi}{2}})^{2n+3}  }  \end{equation}\\

\begin{equation}  A_{4}=  \sum_{n=1}^{\infty}{ E_{n}(4) ({\frac{\pi}{2}})^{2n+3}  }    \end{equation}
where  $  E_{n}(4) =   \frac{D_{n}(4)}{2} + \frac{E_{n}(3)}{!{1}} - \frac{E_{n}(1)}{!{3}} $ \\

$\textbf{{Step 5}}$   \\ 
Integrating Eq.14  over  $ [0,  \theta]$  in Step 5, 

%\begin{equation}   -[\frac{\cos(\theta)}{1^{5}}-\frac{\cos(2\theta)}{2^{5}}+ \frac{\cos(3\theta)}{3^{5}}-\frac{\cos(4\theta)}{4^{5}}+...]+[\frac{1}{1^{5}}-\frac{1}{2^{5}}+\frac{1}{3^{5}}-\frac{1}{4^{5}}+...] =   
%A_{3}\frac{\theta^{2}}{!{2}} - A_{1}\frac{\theta^{4}}{!{4}}  + \frac{1}{2} \sum_{n=1}^{\infty}{ D_{n}(5) (\theta)^{2n+4}  }  
%\end{equation}

\begin{eqnarray*}
  -[\frac{\cos(\theta)}{1^{5}}-\frac{\cos(2\theta)}{2^{5}}+ \frac{\cos(3\theta)}{3^{5}}-\frac{\cos(4\theta)}{4^{5}}+...]+[\frac{1}{1^{5}}-\frac{1}{2^{5}}+\frac{1}{3^{5}}-\frac{1}{4^{5}}+...] =  \\
 A_{3}\frac{\theta^{2}}{!{2}} - A_{1}\frac{\theta^{4}}{!{4}}  + \frac{1}{2} \sum_{n=1}^{\infty}{ D_{n}(5) (\theta)^{2n+4}  }  
\end{eqnarray*}

\begin{equation} \end{equation}

Putting  $ \theta = \frac{\pi}{2} $, we  get,  \\
\begin{equation}  A_{5}= [\frac{1}{1^{5}}-\frac{1}{2^{5}}+\frac{1}{3^{5}}-\frac{1}{4^{5}}...]  = \frac{1}{1-\frac{1}{2^{5}}}[\frac{A_{3}}{!{2}} (\frac{\pi}{2})^{2} - \frac{A_{1}}{!{4}} (\frac{\pi}{2})^{4} + \frac{1}{2} \sum_{n=1}^{\infty}{ D_{n}(5) ({\frac{\pi}{2}})^{2n+4}  } ] \end{equation}\\

\begin{equation}  A_{5}=  \sum_{n=1}^{\infty}{ E_{n}(5) ({\frac{\pi}{2}})^{2n+4}  }    \end{equation}
where  $  E_{n}(5) = \frac{ \frac{D_{n}(5)}{2} + \frac{E_{n}(3)}{!{2}} -  \frac{E_{n}(1)}{!{4}} }{1-\frac{1}{2^{5}}}$ \\

In general, we can derive the following results Eq.20 to Eq.23 using the principle of mathematical induction as shown in Section 2. In Section 3, it is shown that the series
expansion of $A_{2k}$ and $A_{2k+1}$ converges. For $k=1,2,3,...$\\

\begin{equation}  A_{2k}=   [\frac{1}{1^{2k}}-\frac{1}{3^{2k}}+\frac{1}{5^{2k}}+...]  = \sum_{n=1}^{\infty}{ E_{n}(2k) ({\frac{\pi}{2}})^{2n+2k-1}  }    \end{equation}
\begin{equation}  A_{2k+1}=   [\frac{1}{1^{2k+1}}-\frac{1}{2^{2k+1}}+\frac{1}{3^{2k+1}}-\frac{1}{4^{2k+1}}...]  =\sum_{n=1}^{\infty}{ E_{n}(2k+1) ({\frac{\pi}{2}})^{2n+2k}  }    \end{equation}

\begin{equation}  E_{n}(2k)= \frac{(-1)^{k}D_{n}(2k)}{2} + (-1)^{k+1} \sum_{r=0}^{k-1}{ \frac{ (-1)^{r} E_{n}(2r+1)}{!({2(k-r)-1})} }    \end{equation}

\begin{equation}  E_{n}(2k+1)=\frac{ \frac{(-1)^{k}D_{n}(2k+1)}{2} + (-1)^{k+1} \sum_{r=0}^{k-1}{ \frac{ (-1)^{r} E_{n}(2r+1)}{!({2(k-r)})} }  }{1-\frac{1}{2^{2k+1}}}  \end{equation}

Thus we see that $A_{2k}$ and $A_{2k+1}$ can be expressed as an infinite summation of integer powers of $\pi$ with rational coefficients $E_{n}(2k)$ and $E_{n}(2k+1)$.  And, we have \\
\begin{equation}  B_{2k+1}=   [\frac{1}{1^{2k+1}}+\frac{1}{2^{2k+1}}+\frac{1}{3^{2k+1}}+\frac{1}{4^{2k+1}}...]  = \frac{A_{2k+1} }{1-2^{-2k}}   \end{equation}

Using the above results, we can deduce specific values of \textbf{Apery's constant $\zeta{(3)}$ and 
Catalan's constant K} as follows:

\begin{equation}  K = \sum_{n=0}^{\infty} \frac{(-1)^{n}}{(2n+1)^{2}} =   A_{2} =  \sum_{n=1}^{\infty}{ E_{n}(2) ({\frac{\pi}{2}})^{2n+1}  }  \end{equation}

\begin{equation}  \zeta{(3)} = \sum_{n=1}^{\infty} \frac{1}{{n}^{3}} =   B_{3} =  \frac{\sum_{n=1}^{\infty}{ E_{n}(3) ({\frac{\pi}{2}})^{2n+2}  }}{1-2^{-2}}  \end{equation}

\section{\label{sec:level1}Section 2\protect\\  \lowercase{} }

Let us assume that $A_{2k}, A_{2k+1}, E_n(2k), E_n(2k+1)$ are given by Eq. 20 - Eq.23. for some k and that Eq.14 and Eq.17 can be generalized as follows(Inductive hypothesis).
 Let $S_{1}(k)= [\frac{\sin(\theta)}{1^{2k}}-\frac{\sin(2\theta)}{2^{2k}}+ \frac{\sin(3\theta)}{3^{2k}}-\frac{\sin(4\theta)}{4^{2k}}+...]$ and let $S_{2}(k)=  [\frac{\cos(\theta)}{1^{2k+1}}-\frac{\cos(2\theta)}{2^{2k+1}}+ \frac{\cos(3\theta)}{3^{2k+1}}-\frac{\cos(4\theta)}{4^{2k+1}}+...]$\\

\begin{equation} S_{1}(k)=  (-1)^{k} [ \frac{1}{2} \sum_{n=1}^{\infty}{ D_{n}(2k) \theta^{2n+2k-1}  }  ] +  \sum_{r=0}^{k-1}{ \frac{ (-1)^{k-r-1} A_{2r+1} \theta^{2k-2r-1} }{!({2k-2r-1})} }   \end{equation}

\begin{equation} S_{2}(k)=  (-1)^{k+1} [ \frac{1}{2} \sum_{n=1}^{\infty}{ D_{n}(2k+1) \theta^{2n+2k}  }  ] +  \sum_{r=0}^{k}{ \frac{ (-1)^{k-r} A_{2r+1} \theta^{2k-2r} }{!({2k-2r})} }   \end{equation}

We will prove that above equations hold true for k=k+1(Inductive result).
Integrating above equation 28 once from $ [0, \theta] $, we get
$S_{3}(k)=  [\frac{\sin(\theta)}{1^{2k+2}}-\frac{\sin(2\theta)}{2^{2k+2}}+ \frac{\sin(3\theta)}{3^{2k+2}}-\frac{\sin(4\theta)}{4^{2k+2}}+...] $

\begin{equation} S_{3}(k)=   (-1)^{k+1} [ \frac{1}{2} \sum_{n=1}^{\infty} \frac{ D_{n}(2k+1) \theta^{2n+2k+1}  } {2n+2k+1} ] +  \sum_{r=0}^{k}{ \frac{ (-1)^{k-r} A_{2r+1} \theta^{2k-2r+1} }{(2k-2r+1) [!({2k-2r})]} }   \end{equation}

Since $\frac{ D_{n}(2k+1)}{2n+2k+1} = D_{n}(2k+2) $ and $ (2k-2r+1) [!({2k-2r})] = !({2k-2r+1})  $, 
above equation $S_{3}(k)= S_{1}(k+1) $, thus proving the inductive hypothesis for k=k+1. \\

Similarly, Integrating  $ S_{3}(k) $ once from $ [0, \theta] $, we get

\begin{equation} F_{4}(k)= -[\frac{\cos(\theta)}{1^{2k+3}}-\frac{\cos(2\theta)}{2^{2k+3}}+ \frac{\cos(3\theta)}{3^{2k+3}}-\frac{\cos(4\theta)}{4^{2k+3}}+...]+[ \frac{1}{1^{2k+3}}-\frac{1}{2^{2k+3}}+ \frac{1}{3^{2k+3}}-\frac{1}{4^{2k+3}}+...] \end{equation}
\begin{equation} F_{4}(k)= (-1)^{k+1} [ \frac{1}{2} \sum_{n=1}^{\infty}{ D_{n}(2k+2) \frac{\theta^{2n+2k+2}}{2n+2k+2}  }  ] +  \sum_{r=0}^{k}{ \frac{ (-1)^{k-r} A_{2r+1} \theta^{2k-2r+2} }{(2k-2r+2) [!({2k-2r+1})] } }   \end{equation}

Since $\frac{ D_{n}(2k+2)}{2n+2k+2} = D_{n}(2k+3) $ and $ (2k-2r+2) [!({2k-2r+1})] = !({2k-2r+2})  $, and
$ [ \frac{1}{1^{2k+3}}-\frac{1}{2^{2k+3}}+ \frac{1}{3^{2k+3}}-\frac{1}{4^{2k+3}}+...]  = A_{2k+3}$
we can define $S_{4}(k)= A_{2k+3}-F_{4}(k)$ as follows. \\

\begin{equation} S_{4}(k)= A_{2k+3}-F_{4}(k)= [\frac{\cos(\theta)}{1^{2k+3}}-\frac{\cos(2\theta)}{2^{2k+3}}+ \frac{\cos(3\theta)}{3^{2k+3}}-\frac{\cos(4\theta)}{4^{2k+3}}+...] \end{equation}
\begin{equation} S_{4}(k)= A_{2k+3} + (-1)^{k+2} [ \frac{1}{2} \sum_{n=1}^{\infty}{ D_{n}(2k+3) \theta^{2n+2k+2}  }  ] - \sum_{r=0}^{k}{ \frac{ (-1)^{k-r} A_{2r+1} \theta^{2k-2r+2} }{ [!({2k-2r+2})] } }   \end{equation}

Putting $k=k+1$ in Eq. 28, we have
 
\begin{equation} S_{2}(k+1)=   (-1)^{k+2} [ \frac{1}{2} \sum_{n=1}^{\infty}{ D_{n}(2k+3) \theta^{2n+2k+2}  }  ] - \sum_{r=0}^{k}{ \frac{ (-1)^{k-r} A_{2r+1} \theta^{2k-2r+2} }{!({2k-2r+2})} }  
+ A_{2k+3}
 \end{equation}

We can see that $S_{4}(k)=  S_{2}(k+1) $, thus proving the inductive hypothesis for k=k+1. \\

Substituting $ \theta = \pi/2 $ in above equations 27 and 28, we get

\begin{equation} A_{2k}= [\frac{1}{1^{2k}}- \frac{1}{3^{2k}}+  \frac{1}{5^{2k}}...] = (-1)^{k} [ \frac{1}{2} \sum_{n=1}^{\infty}{ D_{n}(2k) (\frac{\pi}{2})^{2n+2k-1}  }  ] +  \sum_{r=0}^{k-1}{ \frac{ (-1)^{k-r-1} A_{2r+1} (\frac{\pi}{2})^{2k-2r-1} }{!({2k-2r-1})} }   \end{equation}

\begin{equation} A_{2k+1}=  2^{2k+1} [ (-1)^{k+1} [ \frac{1}{2} \sum_{n=1}^{\infty}{ D_{n}(2k+1) (\frac{\pi}{2})^{2n+2k}}  ] +  \sum_{r=0}^{k}{ \frac{ (-1)^{k-r} A_{2r+1} (\frac{\pi}{2})^{2k-2r} }{ !({2k-2r})} } ]  \end{equation}

where $A_{2k+1}= [ \frac{1}{1^{2k+1}}-\frac{1}{2^{2k+1}}+ \frac{1}{3^{2k+1}}-\frac{1}{4^{2k+1}}+...]$. These show that $A_{2k}$ and $A_{2k+1}$ can be expressed as linear combination of powers of $\pi/2$. Let us assume the following inductive hypothesis for some k:

\begin{equation}  A_{2k+1} =\sum_{n=1}^{\infty}{ E_{n}(2k+1) ({\frac{\pi}{2}})^{2n+2k}  }    \end{equation}

\begin{equation}  E_{n}(2k+1)=\frac{ \frac{(-1)^{k}D_{n}(2k+1)}{2} + (-1)^{k+1} \sum_{r=0}^{k-1}{ \frac{ (-1)^{r} E_{n}(2r+1)}{!({2(k-r)})} }  }{1-\frac{1}{2^{2k+1}}}  \end{equation}

We will prove  that this hypothesis holds for k=k+1. Putting k=k+1 in Eq.36,

\begin{equation} A_{2k+3}=  2^{2k+3} [ (-1)^{k} [ \frac{1}{2} \sum_{n=1}^{\infty}{ D_{n}(2k+3) (\frac{\pi}{2})^{2n+2k+2}}  ] +  \sum_{r=0}^{k}{ \frac{ (-1)^{k+1-r} A_{2r+1} (\frac{\pi}{2})^{2k-2r+2} }{ !({2k-2r+2})} } + A_{2k+3} ]  \end{equation}

This can be written as follows:

\begin{equation} A_{2k+3}= \frac{1}{ 1- \frac{1}{2^{2k+3}} } [ (-1)^{k+1} [ \frac{1}{2} \sum_{n=1}^{\infty}{ D_{n}(2k+3) (\frac{\pi}{2})^{2n+2k+2}}  ] -  \sum_{r=0}^{k}{ \frac{ (-1)^{k+1-r} A_{2r+1} (\frac{\pi}{2})^{2k-2r+2} }{ !({2k-2r+2})} }  ]  \end{equation}

Using Eq. 37 to replace $A_{2r+1}$ in above equation, interchanging order of summation and taking out common factor $ (\frac{\pi}{2})^{2n+2k+2}$, we get

\begin{equation}  A_{2k+3} =\sum_{n=1}^{\infty}{ E_{n}(2k+3) ({\frac{\pi}{2}})^{2n+2k+2}  }    \end{equation}

\begin{equation}  E_{n}(2k+3)=\frac{ \frac{(-1)^{k+1}D_{n}(2k+3)}{2} + (-1)^{k} \sum_{r=0}^{k}{ \frac{ (-1)^{r} E_{n}(2r+1)}{!({2(k-r+1)})} }  }{1-\frac{1}{2^{2k+3}}}  \end{equation}

Thus we have proved the inductive result for k=k+1. Eq. 37 and Eq. 38 imply above results.\\

Similarly, $A_{2k}$ in  Eq.35 can be extended to  $A_{2k+2}$ by replacing k with k+1:

\begin{equation} A_{2k+2} = (-1)^{k+1} [ \frac{1}{2} \sum_{n=1}^{\infty}{ D_{n}(2k+2) (\frac{\pi}{2})^{2n+2k+1}  }  ] +  \sum_{r=0}^{k}{ \frac{ (-1)^{k-r} A_{2r+1} (\frac{\pi}{2})^{2k-2r+1} }{!({2k-2r+1})} }   \end{equation}

 Let us assume the following inductive hypothesis for $A_{2k}$ and $E_{n}(2k)$  some k:

\begin{equation}  A_{2k}=   [\frac{1}{1^{2k}}-\frac{1}{3^{2k}}+\frac{1}{5^{2k}}+...]  = \sum_{n=1}^{\infty}{ E_{n}(2k) ({\frac{\pi}{2}})^{2n+2k-1}  }    \end{equation}

\begin{equation}  E_{n}(2k)= \frac{(-1)^{k}D_{n}(2k)}{2} + (-1)^{k+1} \sum_{r=0}^{k-1}{ \frac{ (-1)^{r} E_{n}(2r+1)}{!({2(k-r)-1})} }    \end{equation}

Using Eq. 37 to replace $A_{2r+1}$ in above equation 43, interchanging order of summation and taking out common factor $ (\frac{\pi}{2})^{2n+2k+1}$, we get

\begin{equation}  A_{2k+2} = \sum_{n=1}^{\infty}{ E_{n}(2k+2) ({\frac{\pi}{2}})^{2n+2k+1}  }    \end{equation}

\begin{equation}  E_{n}(2k+2)= \frac{(-1)^{k+1}D_{n}(2k+2)}{2} + (-1)^{k} \sum_{r=0}^{k}{ \frac{ (-1)^{r} E_{n}(2r+1)}{!({2(k-r)+1})} }    \end{equation}

Thus we have proved the inductive result for k=k+1. Eq. 44 and Eq. 45 imply above results.\\

\section{\label{sec:level1}Section 3\protect\\  \lowercase{} }

In this section, it will be shown that the series expansion of $A_{2k}$ and $A_{2k+1}$ in Eq. 20 and Eq. 21 converges.\\

We know that $D_{n}(k) = \frac{c_{n}}{2^{2n-1}(2n+k-1)P_{k}}$, which can be rewritten as follows:

\begin{equation} D_{n}(k) = \frac{D_{n}(1)}{(2n+k-1)P_{k-1}}  \end{equation}

We will write $E_{n}(2k)$ and $E_{n}(2k+1)$ in Eq. 22 and Eq.23 in terms
of $D_{n}(1)$ as follows:

\begin{equation} E_{n}(1) = D_{n}(1) \end{equation}

\begin{equation} E_{n}(2) = E_{n}(1) - \frac{D_{n}(2)}{2} =  D_{n}(1) F_{n}(2) \end{equation}
where $F_{n}(2) = 1 - \frac{1}{2 (2n+1)}$ and $\lim_{n\to\infty} {F_{n}(2) = K(2)}$ where $K(2)$ is a constant.\\

\begin{equation} E_{n}(3) = \frac{ - \frac{D_{n}(3)}{2} + \frac{E_{n}(1)}{!{2}}}{1-\frac{1}{2^{3}}} = D_{n}(1) F_{n}(3) \end{equation}
where $F_{n}(3) = \frac{\frac{1}{!{2}} - \frac{1}{2 ((2n+2) P_{2})}}{1-\frac{1}{2^{3}}}$ and $\lim_{n\to\infty} {F_{n}(3) = K(3)}$ where $K(3)$ is a constant.\\

\begin{equation} E_{n}(4) =   \frac{D_{n}(4)}{2} + \frac{E_{n}(3)}{!{1}} - \frac{E_{n}(1)}{!{3}} = D_{n}(1) F_{n}(4) \end{equation}

where $F_{n}(4) = -\frac{1}{!{3}}+  \frac{1}{2 ((2n+3) P_{3})} + \frac{\frac{1}{!{2}} - \frac{1}{2 ((2n+2) P_{2})}}{1-\frac{1}{2^{3}}}  $ and $\lim_{n\to\infty} {F_{n}(4) = K(4)}$ where $K(4)$ is a constant.\\

\begin{equation} E_{n}(5) = \frac{ \frac{D_{n}(5)}{2} + \frac{E_{n}(3)}{!{2}} -  \frac{E_{n}(1)}{!{4}} }{1-\frac{1}{2^{5}}} =  D_{n}(1) F_{n}(5) \end{equation}

where $F_{n}(5) = \frac{-\frac{1}{!{4}} + \frac{\frac{1}{!{2}} - \frac{1}{2 ((2n+2) P_{2})}}{(!{2})(1-\frac{1}{2^{3}})} + \frac{1}{2 ((2n+4) P_{4})}}{1-\frac{1}{2^{5}}}$
and $\lim_{n\to\infty} {F_{n}(5) = K(5)}$ where $K(5)$ is a constant.\\

Let us assume the \textbf{Inductive Hypothesis that $E_{n}(2k) = D_{n}(1) F_{n}(2k)$
and $E_{n}(2k+1) = D_{n}(1) F_{n}(2k+1)$ } and that $\lim_{n\to\infty} {F_{n}(2k) = K(2k)}$ and $\lim_{n\to\infty} {F_{n}(2k+1) = K(2k+1)}$ where $K(2k)$ and 
$K(2k+1)$ are constants. Substituting this in Eq.22 and Eq.23, we can write
Eq.42 and Eq.47 as $E_{n}(2k+2) = D_{n}(1) F_{n}(2k+2)$
and $E_{n}(2k+3) = D_{n}(1) F_{n}(2k+3)$ where $\lim_{n\to\infty} {F_{n}(2k+2) = K(2k+2)}$ and $\lim_{n\to\infty} {F_{n}(2k+3) = K(2k+3)}$ where $K(2k+2)$ and 
$K(2k+3)$ are constants, thus proving the Inductive Result.\\

Hence we can write
\begin{equation} \lim_{n\to\infty} \frac{E_{n+1}(2k)}{E_{n}(2k)} = K_{0} \lim_{n\to\infty} \frac{D_{n+1}(1)}{D_{n}(1)} = \frac{K_{0}}{4} \lim_{n\to\infty} \frac{c_{n+1}}{c_{n}} \frac{2n}{2n+2} \end{equation}

where $K_{0} =  \lim_{n\to\infty} \frac{F_{n+1}(2k)}{F_{n}(2k)} = 1$. 
Given that $\lim_{n\to\infty} |\frac{c_{n+1}}{c_{n}}| < 1 $ in the series expansion of 
$\tan(\theta)$ in Eq. 3, we get the result\\

\begin{equation} \lim_{n\to\infty} |\frac{E_{n+1}(2k)}{E_{n}(2k)}| < 1   \end{equation}

Similarly it can be shown that 
\begin{equation} \lim_{n\to\infty} |\frac{E_{n+1}(2k+1)}{E_{n}(2k+1)}| < 1 \end{equation}

Hence the the series expansion of $A_{2k}$ and $A_{2k+1}$ in Eq. 20 and Eq. 21 converges.

\section{\label{sec:level1} Conclusion \protect\\  \lowercase{} }

It has been shown that Riemann's Zeta function with odd positive integer argument 
can be represented as an infinite summation of integer powers of $\pi$ with rational coefficients. Specific
values for Apery's Constant and Catalan's Constant have been derived.

\section{\label{sec:level1}References\protect\\  \lowercase{} }

[1] Hardy, G. H. Divergent Series. New York: Oxford University Press, 1949. 

[2] Abramowitz, M. and Stegun, I. A. (Eds.). "Circular Functions." §4.3 in Handbook of Mathematical Functions with Formulas, Graphs, and Mathematical Tables, 9th printing. New York: Dover, p. 75, 1972. 


\bibliographystyle{elsarticle-num}
\bibliography{<your-bib-database>}






% \end{thebibliography}

\end{document}